\documentclass[11pt,a4paper]{article}
\usepackage{charter}                                        %
\usepackage{amsfonts,amsmath,amssymb}                     %
\usepackage[a4paper]{geometry}                            %
\usepackage{fancyhdr}                                     %
\usepackage{color}                                        %
\usepackage[pdftex]{hyperref} 
\usepackage{graphicx}
\hypersetup{                                              %
    a4paper,                                              %
    breaklinks                                            %
}                                                         %
{\ \rule{0.5em}{0.5em}}                                   %
\usepackage{tabularx}
\usepackage{booktabs}
\usepackage{quoting}
\usepackage[all]{nowidow}
\begin{document}
%
%
\title{\vspace*{-1cm} Flipped classroom in Introductory Linear Algebra by utilizing\\ Computer Algebra System {\sl SageMath} and a free electronic book}%
%
%
\author{\begin{tabular}{c}
{\normalsize \textit{N. Karjanto and S.-G. Lee}} \\
{\small Department of Mathematics, University College} \\
{\small Natural Science Campus, Sungkyunkwan University} \\
{\small Suwon 16419, } 
{\small Republic of Korea} \end{tabular}
}%
%
\date{}                                                   %
\maketitle                                                %
%
%
\begin{abstract}
This article describes the authors' teaching experience in flipping the class of a basic undergraduate mathematics course Introductory Linear Algebra. We utilize a full-featured free electronic textbook, online lecture notes, an intranet learning management system (LMS) {\sl icampus}, the video-sharing website {\sl YouTube} and a Computer Algebra System (CAS) {\sl SageMath} in our flipped classroom approach. Prior coming to the class, the students are assigned to complete a portion of reading assignments from the free e-book and to watch recorded video lectures from {\sl YouTube} that cover the relevant materials. Announcements are made through the LMS {\sl icampus}, as well as disseminating teaching materials, assigning tasks, and fostering online communication outside the class hours. We dedicated in-class activities with problem-solving sessions, question and answer, discussion and providing instant feedback to the students. Hence, both out-of-class and in-class sessions are filled with active learning and interactive engagement. From the students' feedback and opinion, we discovered that some of them preferred the traditional classroom rather than the flipped classroom pedagogy. A possible reason is either they are used to the traditional lecture style where receiving information occurs inside the classroom or the flipped classroom approach that we implement does not satisfy students' learning styles. Nevertheless, several students enjoyed problem-solving activities where they have an opportunity to communicate with the instructors whenever they encounter some difficulty.
\end{abstract}%
%
%
\thispagestyle{fancy}                                     %
%
%
\section{Introduction}
As educators and practitioners, we strive our best in assisting our students in their learning process. This effort will not only help them to retain the knowledge they learn but also to develop their problem-solving and critical thinking skills that will be useful even after they completed a certain course. Adopting some kind of active learning is one such strategy for better instruction and flipping the class is one feature of active learning which has become popular during the past decade. In this article, we will discuss flipped classroom for an introductory course in Linear Algebra by utilizing a full-featured free electronic textbook and CAS {\sl SageMath}.

The term `flipped classroom' is more popular for K-12 education. It was coined and popularized by Chemistry teachers Jonathan Bergmann and Aaron Sams from Woodland Park High School, Colorado back in 2007~\cite{bergmann2012flip,bergmann2014flipped}. At the college level, the term `inverted classroom' is more commonly used. The term was coined by a group of Economics faculty at Miami University, Ohio in 2000~\cite{lage2000inverting}. Although the term might be new, the idea borrows from the centuries-old method of assigning reading materials to students before they come to the class. If students are prepared beforehand, in-class time can be dedicated to active learning and higher-order skills, including problem-solving, discussion and critical thinking.

There are a number of studies on the flipped classroom in the teaching and learning of Mathematics, covering the primary, secondary and tertiary educations. The following provides some examples of the flipped classroom in a number of mathematical subjects: PreCalculus, Algebra and Geometry~\cite{fulton2012upside}, College Algebra~\cite{overmyer2014flipped,van2015adventures,acelajado2017flipped}, Calculus~\cite{mcgivney2013flipping,jungic15,sahin2015flipping,sonnert2015impact,maciejewski2016flipping}, Multivariable Calculus~\cite{ziegelmeier2015flipped}, Statistics~\cite{wilson2013flipped} and Actuarial Science~\cite{butt2014student}.

There exists research of flipped classroom in Linear Algebra and the list is not exhaustive.
Several flipped classroom designs which include flipping a single-topic, an entire course and in the form of a series of the workshop has been proposed~\cite{talbert2014inverting}.
Literature has also provided an evidence of the flipped classroom effectiveness in Linear Algebra.
For instance, a comparison between the flipped classroom and the traditional lecture can be found in~\cite{love2014student},
where students' understanding is better in the former.
Another study on students' academic performance and perception is given in~\cite{murphy2016student},
where students in the flipped classroom achieved superior comprehension and exhibit a positive attitude in terms of enjoyment and confidence. 
See also~\cite{novak2017flip} and~\cite{parklee2016} for further supporting evidence of positive students' perspective, participation, interest and self-directed learning skills in the flipped classroom environment.

In this article, we share our teaching experience in flipping the class for an introductory course in Linear Algebra and attempt to fill the gap of flipped classroom study in this area. There are a number of factors that make this study is unique and interesting, i.e. utilization of a free electronic textbook, an adaptation of CAS {\sl SageMath} and incorporation of the learning management system (LMS) {\sl icampus} in the teaching and learning. In addition, the class environment is an English-medium instruction where both the instructors and the majority of our students are not native English speakers. 

The organization of this article is given as follows. 
After this Introduction, Section~\ref{LA} explains the course logistics, the LMS \textsl{i-campus}, the free electronic textbook and online lecture notes.
Section~\ref{flcl} provides features of the flipped classroom in Introductory Linear Algebra course.
This includes video recordings, CAS {\sl SageMath} and some flipped classroom activities.
Section~\ref{feedback} describes the students' feedback and perception in connection with the flipped classroom. A reflection from instructors is also discussed. Finally, Section~\ref{conclusion} provides conclusion and remark to our findings.

\section{Introductory Linear Algebra} \label{LA}

\subsection{Course logistics}

Sungkyunkwan University (SKKU) is a private university in the Republic of Korea (South Korea) affiliated with Samsung Corporation with its history dated back to 1398.
It has two campuses: the Natural Science Campus (NSC) and the Humanities and Social Science Campus (HSCC), located in Suwon and Seoul, respectively.
Linear Algebra (course code GEDB003) is a three-credit introductory course offered as a part of Basic Science and Mathematics modules at SKKU.
It is a service course offered by Department of Mathematics at the NSC for students from various majors in Science, Computer Science and Engineering. Depending on the demand, a typical number of sections offered each semester ranges from four to six, with up to two sections being offered in the HSCC campus.
One instructor acts as a course coordinator and organizes for possible common examinations.

Three out of six sections offered during Spring 2016, three sections are offered as the flipped classroom, whereas the other three are delivered in the traditional lecture style. The classes meet two times a week with each meeting lasts for 75~minutes. The number of students per class is limited to 70.
The classroom is a typical lecture hall with a fixed seating arrangement with a camera installed at the middle or at the rear of the classroom.
Instructors may choose to record in-class sessions and the recordings will be uploaded automatically into the university's database system.
More detailed information on each section of this class during Spring 2016 is given in Table~\ref{section}.

\subsection{LMS {\sl icampus}} 

\begin{figure}[!ht]
\begin{center}
\includegraphics[width = 0.45\textwidth]{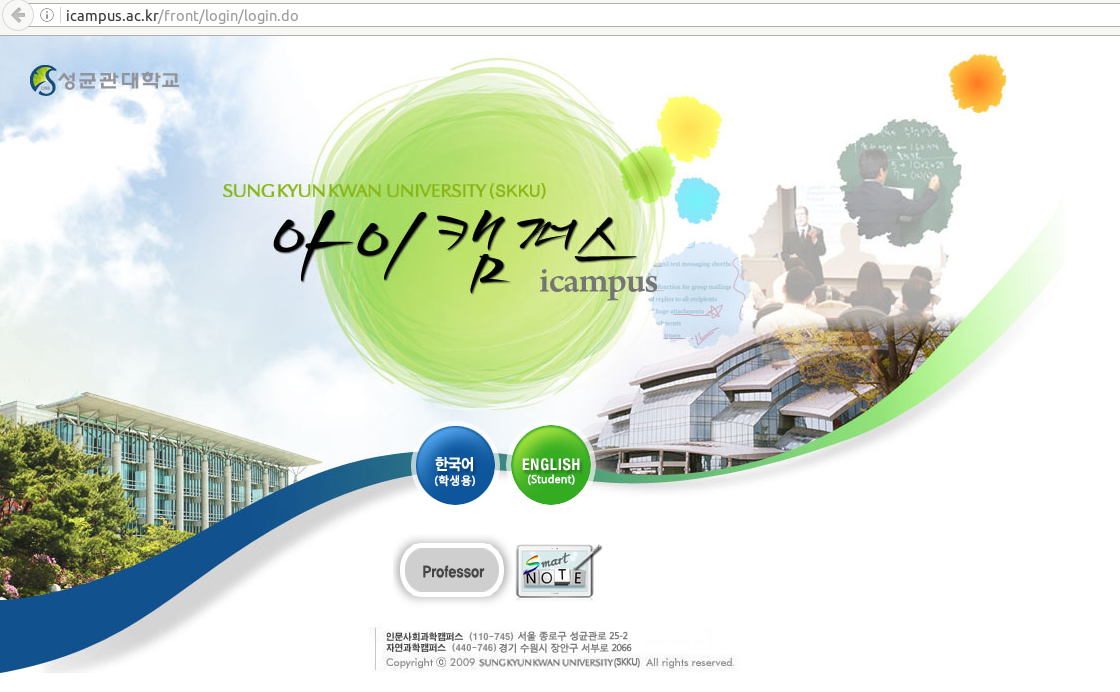} \hspace{1cm}
\includegraphics[width = 0.45\textwidth]{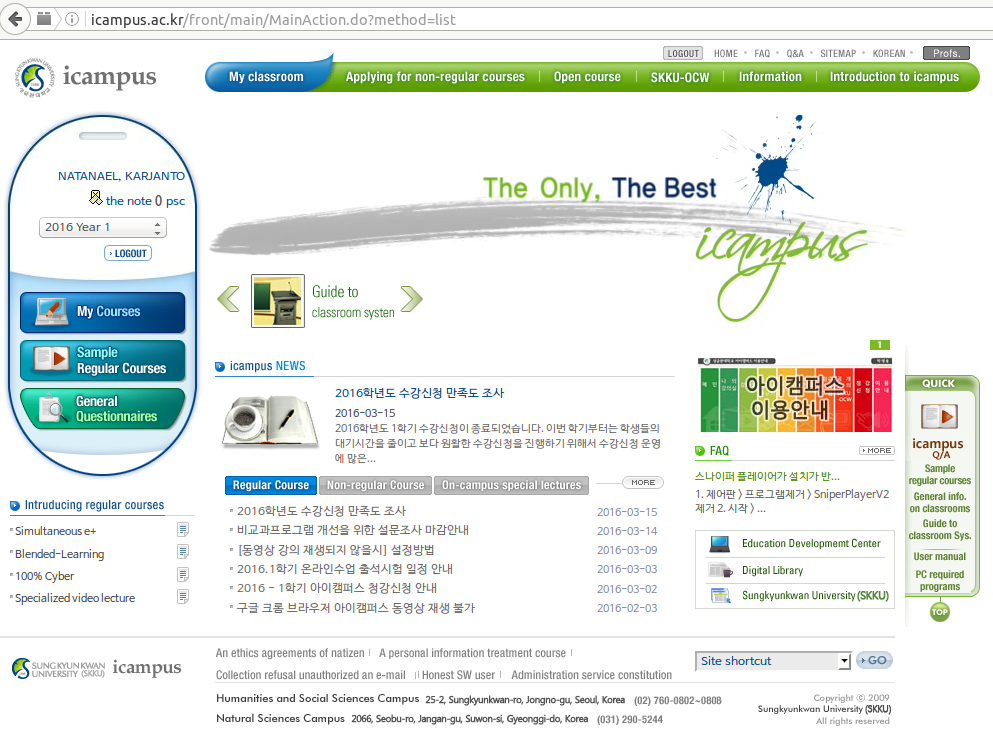} \vspace{0.5cm}\\
\includegraphics[width = 0.45\textwidth]{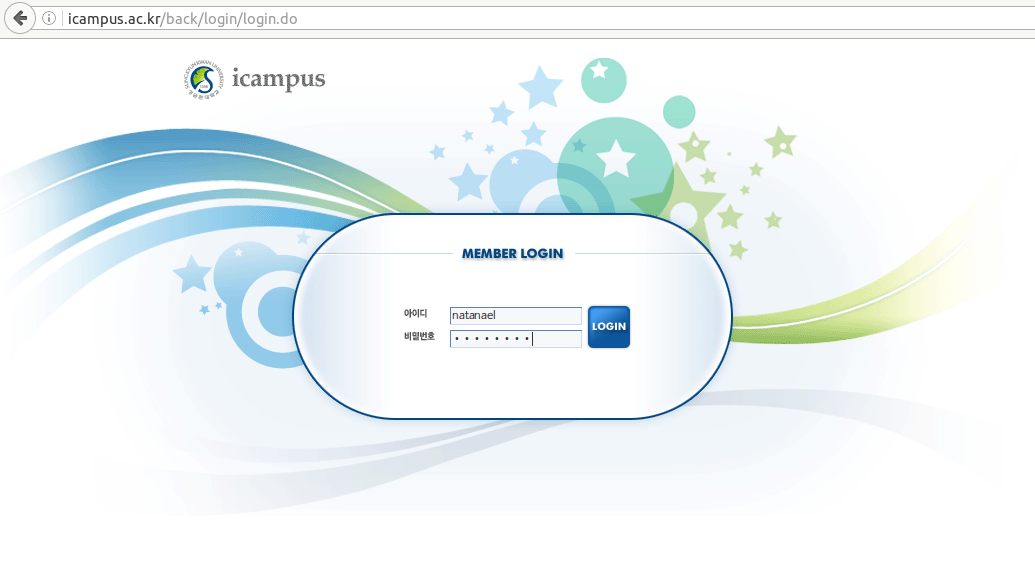} \hspace{1cm}
\includegraphics[width = 0.45\textwidth]{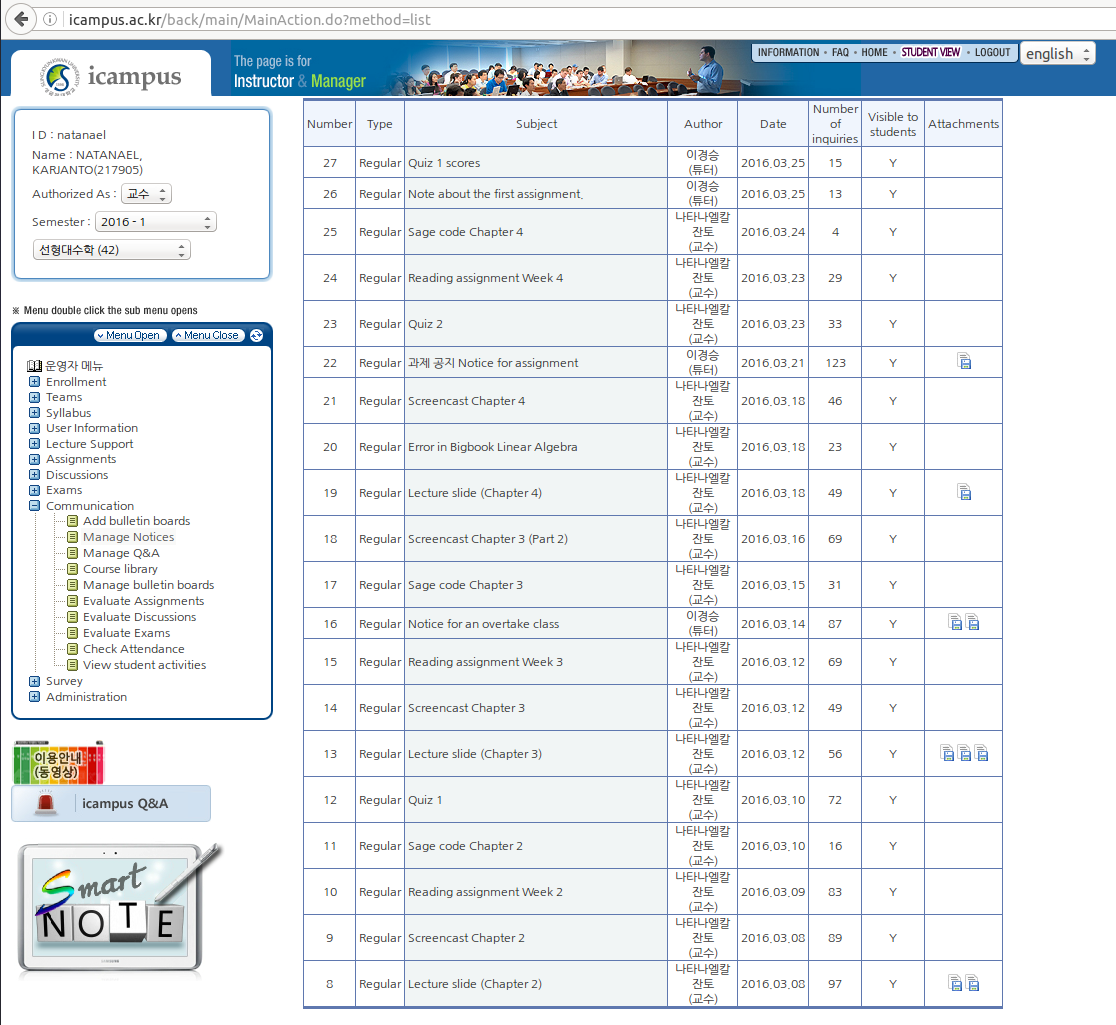}
\end{center}
\caption{(Top panels) Screenshots of the LMS {\sl icampus} homepage, before login (left panel) and after login (right panel).
(Bottom panels) Screenshots of the LMS {\sl icampus} page for faculty, instructors and managers, before login (left panel) and after login (right panel).} \label{icampus}
\end{figure}

Out-of-class communication between instructors and students is conducted using an LMS {\sl icampus}.
This LMS has features which are similar to other well-established LMS, such as {\sl Moodle}, {\sl Edmodo}, or {\sl Blackboard Learn}.
Instructors can write announcements, upload course materials, assign tasks and initiate open discussion. 
Students can do the latter as well as sending messages to their instructor. 
An excellent feature of {\sl icampus} is the possibility to access and to watch video recordings uploaded earlier.
Screenshots of {\sl icampus} before and after the login pages are displayed in Figure~\ref{icampus}.
{\renewcommand{\baselinestretch}{1} 
\begin{table}[h]
\begin{center}
\begin{tabular}{@{}lclllcc@{}}
\toprule
Campus & Section & Instructor    & Day      & Time         & Flipped & Participants    \\ \hline
HSSC   & 01      & Sang Woon Yun & Tuesday  & 12:00--13:15 & No      & 48              \\
       &         &               & Thursday & 15:00--16:15 &         &                 \\ \hline
HSSC   & 02      & Sang Woon Yun & Tuesday  & 15:00--16:15 & No      & 26              \\
       &         &               & Thursday & 12:00--13:15 &         &                 \\ \hline
 NSC   & 41      & SGL           & Tuesday  & 09:00--10:15 & Yes     & 21 		       \\ 
       &         &               & Thursday & 10:30--11:45 &         &				   \\ \hline  
 NSC   & 42      & NK            & Tuesday  & 10:30--11:45 & Yes     & 70 		       \\        
       &         &               & Thursday & 09:00--10:15 &         &   		       \\ \hline
 NSC   & 43      & Lois Simon    & Tuesday  & 12:00--13:15 & No      & 70 		       \\        
       &         &               & Thursday & 15:00--16:15 &         &   		       \\ \hline      
 NSC   & 44      & NK            & Tuesday  & 13:30--14:45 & Yes     & 70 		       \\        
       &         &               & Thursday & 16:30--17:45 &         &   		       \\ 
\bottomrule 
\end{tabular}
\end{center}
\caption{An information on sections of Linear Algebra offered in both campuses during Spring 2016.} \label{section}
\end{table}}

\vspace*{-0.5cm}
\subsection{Free electronic book} \label{ebook}

A new model of interactive digital mathematics textbook has been introduced~\cite{lee2017textbook}. 
The subject covers a number of mathematics courses, including Calculus, Linear Algebra, Differential Equations, Probability and Statistics, and Engineering Mathematics. In particular, they discussed a digital textbook on Linear Algebra and focus on their experience in using digital contents and interactive laboratories for developing a new model for a digital textbook in Linear Algebra. They also include some appraisal from the readers who feel that the digital textbook is very useful for them, a staff employed by the Bank of Korea and a faculty member from University of California at Berkeley in the US.

We adopt a textbook developed by our own faculty and the electronic version of this book is freely available from the web address {\small \url{http://ibook.skku.edu/Viewer/LA-Text-Eng}} and the Portable Document Format (PDF) file can be downloaded from the above address. 
It is accessible to anyone and the students can download and use it for free of charge and thus we are promoting the utilization of electronic book (e-book). The students who wish to obtain a hard copy version of the textbook can order it from {\small \url{http://pod.kyobobook.co.kr}} by `print on demand (POD)' service from Kyobo, the largest bookstore in South Korea, for an affordable price of around 17~USD at the current exchange rate.

\begin{figure}[h]
\begin{center}
\includegraphics[width = 0.25\textwidth]{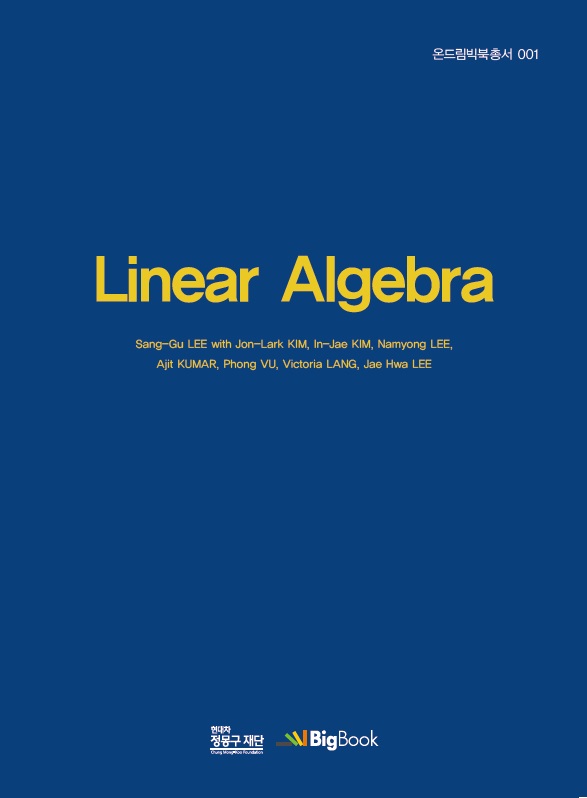} \hspace{2cm}
\includegraphics[width = 0.23\textwidth]{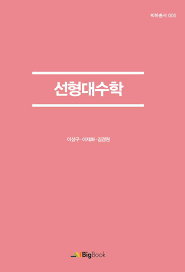}
\end{center}
\caption{The covers of the adopted textbook {\sl Bigbook Linear Algebra}. The corresponding English and Korean versions are shown on the left-hand and the right-hand sides, respectively.}
\label{book}
\end{figure}

The English version of this textbook is based on and is translated from the original version of Korean language. The Korean language version of the e-book is also freely available and accessible online from {\small \url{http://ibook.skku.edu/Viewer/LA-Texbook}}. 
The front cover of both versions is shown in Figure~\ref{book}. This Linear Algebra e-book contains ten chapters, and the title of each chapter is given in the ascending order as follows: Vectors, System of Linear Equations, Matrix and Matrix Algebra, Determinant, Matrix Models, Linear Transformations, Dimension and Subspaces, Diagonalization, General Vector Spaces and Jordan Canonical Form. The appendix of the e-book contains sample exam questions with their solutions and some worked out exercises. These exercises have been solved, revised and finalized by the students and the final check is conducted by the main author of the e-book.

The e-book contains many worked examples that will enhance students' understanding and in learning Linear Algebra. It also contains embedded interactive materials with the CAS {\sl SageMath}. Every section of the book includes links to recorded video lectures on {\sl YouTube} and practice web pages. {\sl SageMath} codes are included in the textbook and thus the students can practice using a computer or a mobile device. They can simply copy and paste the commands and the computer code from the e-book to a {\sl SageMath} worksheet or {\sl SageMath} cell server and modify them accordingly to enhance their understanding of the material. For instance, by playing with a larger size of a matrix or by changing some figures from simple ones to more complicated figures, students are able to perceive the usefulness of the CAS and the aid that it provides.

\vspace*{-0.5cm}
\subsection{Online lecture notes} \label{notes}

In addition to the Linear Algebra textbook in digital format, recently one of us has also developed online lecture notes. These notes are accessible at {\small \url{http://matrix.skku.ac.kr/LA/}}. As can be seen on this website, the lecture notes are arranged according to each chapter of the e-book, and there are also ten chapters in total. Whenever possible, each example is accompanied by {\sl SageMath} code where there is a link to the corresponding web page. The {\sl Chrome} web browser from {\sl Google} is recommended to be used. This is very useful for a class demonstration using {\sl SageMath} as the instructor can simply press the `Evaluate' button and the computational result will be displayed accordingly. This action saves time in comparison to copy and paste the computational codes to a separate cell server.

Furthermore, the {\sl YouTube} video links of the recorded in-class lectures are also displayed below the URL links of the lecture notes.
These video links are arranged according to the section number and related topic systematically. 
This allows the students to access and to watch the recorded video lectures on the relevant topic before they come to the class.
Hence, providing the links which are accessible to the students supports the structure of flipped classroom
where the students can gain basic information outside the classroom by listening to podcasts or watching recorded video lectures.
In-class time can thus be dedicated to more rigorous problem-solving activities, discussion and feedback.
Past sample exams including both midterm and final examinations are also available on the web page.

\vspace*{-0.5cm}
\section{Flipped classroom in Introductory Linear Algebra} \label{flcl}

\subsection{Video recordings}

As mentioned previously in the Introduction, one feature of active learning is adopting and implementing flipped classroom approach in our course. In order to make the flipped classroom successful, the students must prepare prior coming to the class and they must participate actively in in-class activities. For the former, assigned tasks include some reading assignments or watching recorded video lectures. In order to make sure that they come prepared, a small assignment can be embedded into the task and the mark from this assignment is to be included in the final grade calculation.

\begin{figure}[h]
\begin{center}
\includegraphics[width = 0.4\textwidth]{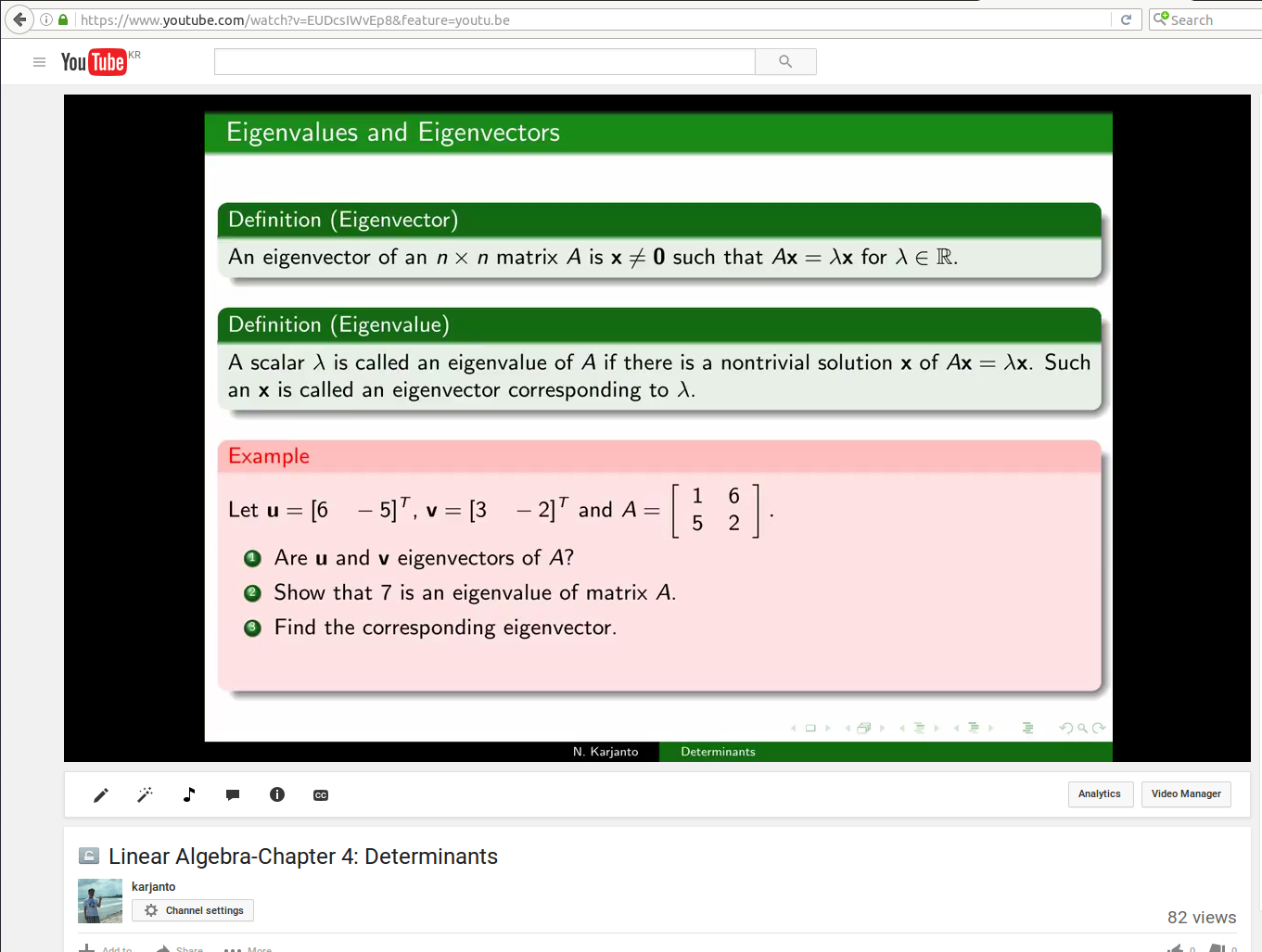} \hspace{1cm}
\includegraphics[width = 0.5\textwidth]{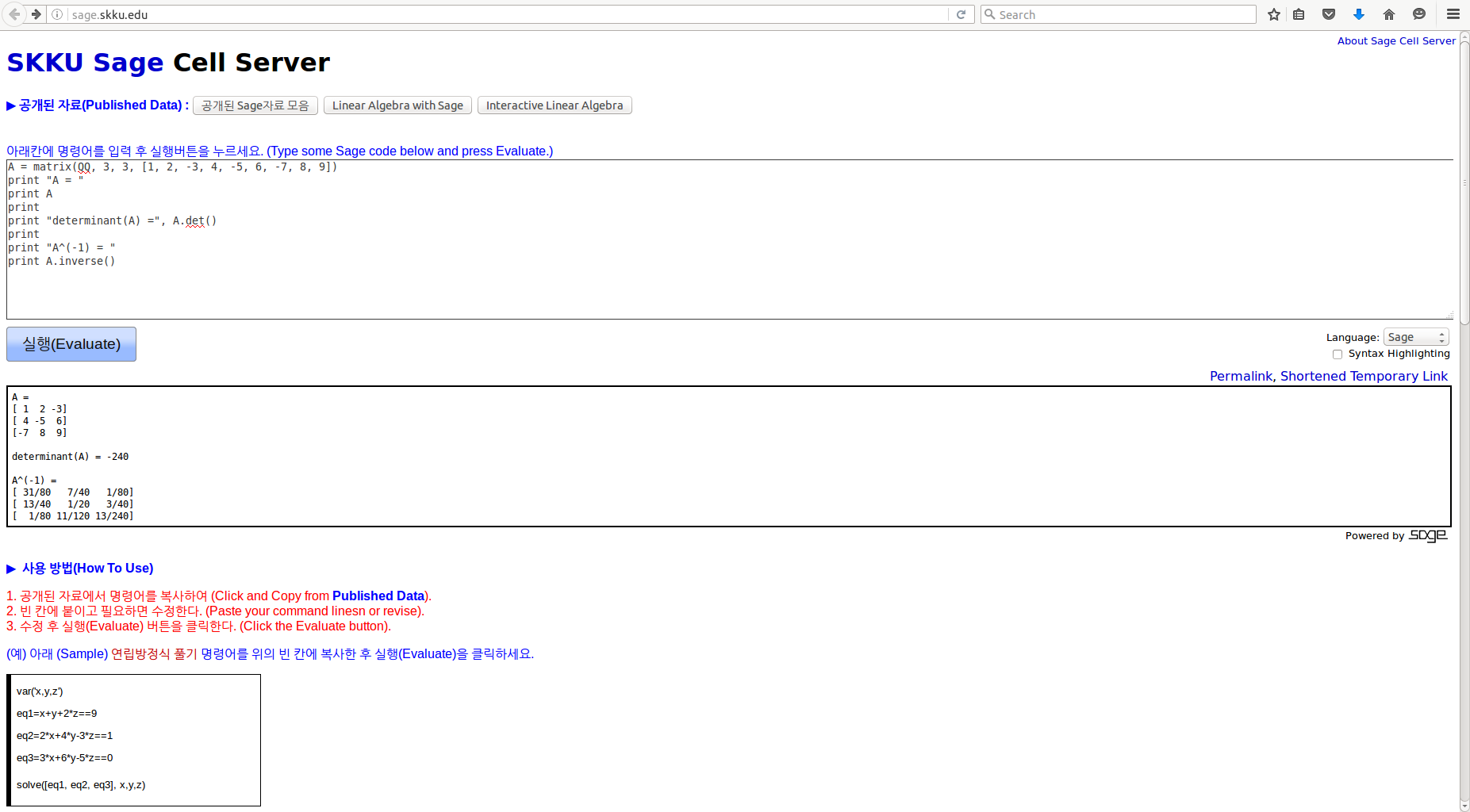}
\end{center}
\caption{(Left) A~screenshot from a {\sl YouTube} page which shows a screencast video material of Linear Algebra on determinant, eigenvalue and eigenvector. (Right) A~screenshot of the SKKU's {\sl Sage} cell server with an example of calculating a determinant of a $3 \times 3$ matrix and the inverse of a matrix.}
\label{sage}
\end{figure}

\begin{figure}[h]
\begin{center}
\includegraphics[width = 0.45\textwidth]{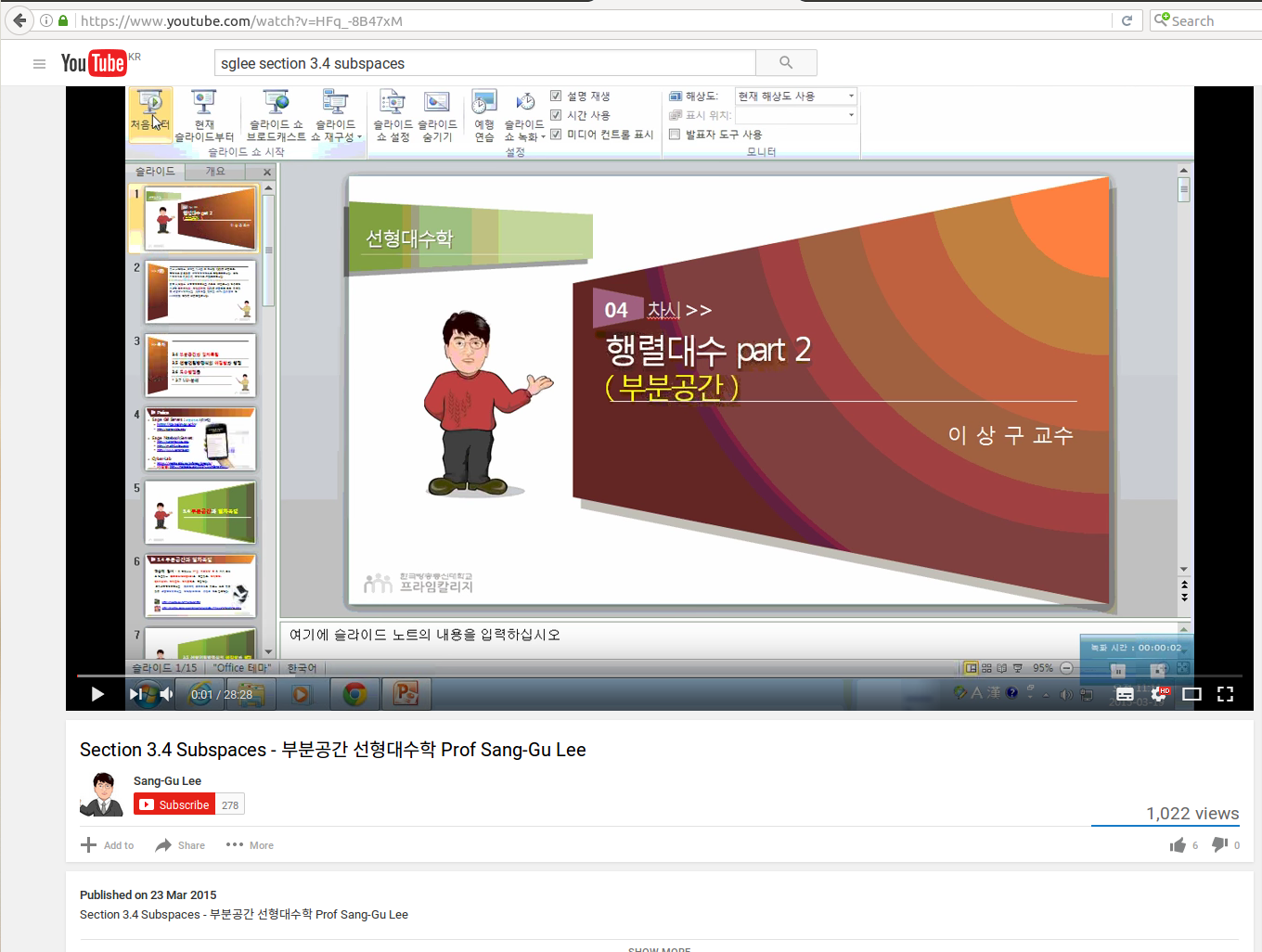} \hspace{1cm}
\includegraphics[width = 0.45\textwidth]{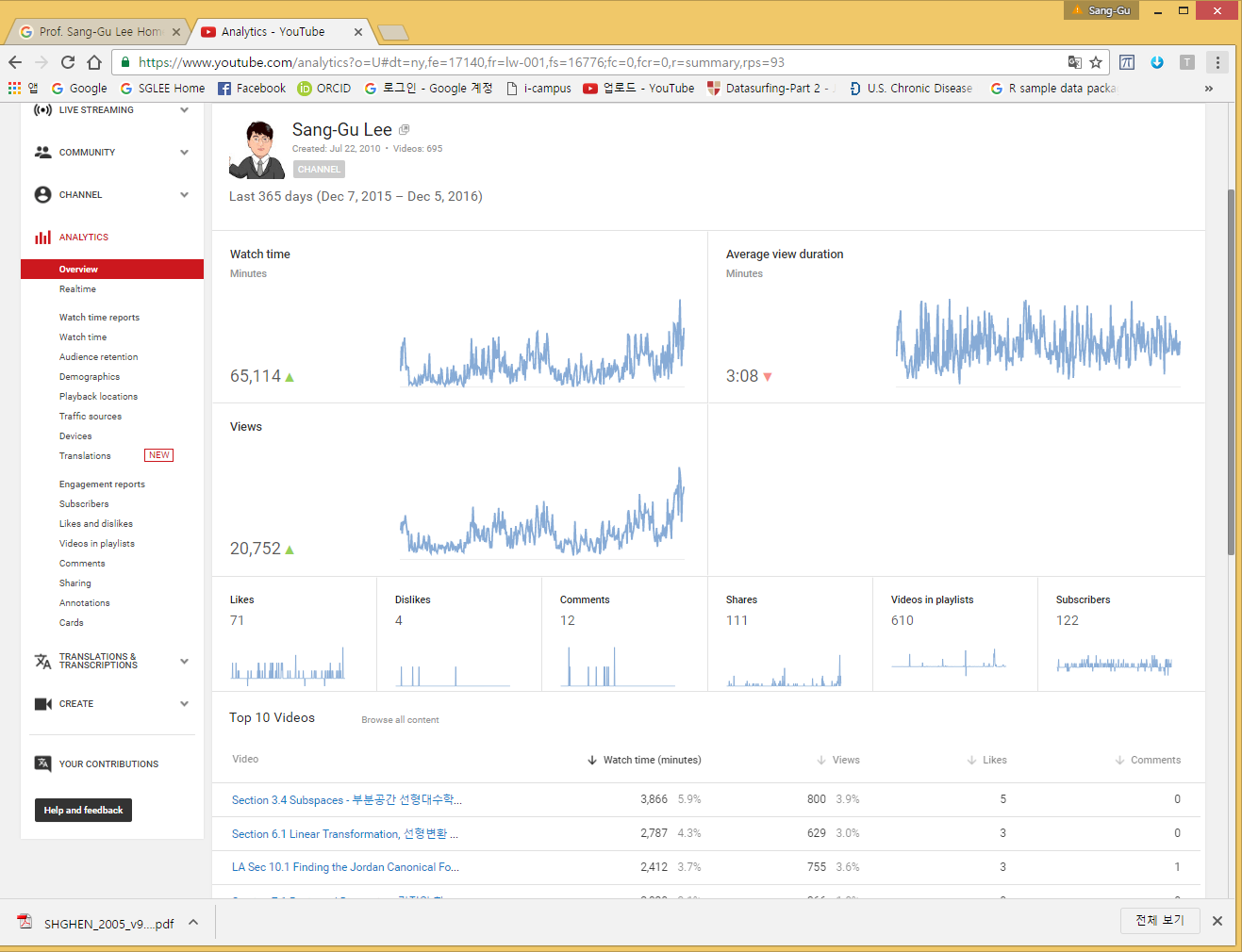}
\end{center}
\caption{(Left) A screenshot from a {\sl YouTube} page which shows a screencast video material of Linear Algebra on subspace. (Right) A screenshot of the Analytics page from {\sl YouTube}.}
\label{analytic}
\end{figure}

When it comes to selecting video lectures to be assigned to the students, we have many options. We can record our own lectures, upload them to a video-sharing website and send the URL link to the students or we can use videos from other instructors available online. We have decided to choose the former. For the latter, one can access them from the video-sharing website {\sl YouTube}. Uploading video recordings online will be useful not only in helping the students to prepare materials but also beneficial in case they are interested to view them again at a later time. In addition, a personalized learning resource {\sl Khan Academy} has more than 2400 free instructional videos on many different subjects. Depending on the materials, the duration of the videos may range from short (around five minutes) to long (around one hour) ones. The screenshot of the video recordings is shown Figures~\ref{sage} and~\ref{analytic}.

The left panel of Figure~\ref{sage} shows a screenshot from a {\sl YouTube} page showing a screencast video material of Linear Algebra on determinant, eigenvalue and eigenvector. The right panel of Figure~\ref{sage} shows a screenshot of the SKKU's {\sl SageMath} cell server with an example of calculating a determinant of a $3 \times 3$ matrix and the inverse of a matrix. The left panel of Figure~\ref{analytic} shows a screenshot from a {\sl YouTube} page which shows a screencast video material of Linear Algebra on subspace. The right panel of Figure~\ref{analytic} displays a screenshot of the Analytics data page from {\sl YouTube}. From the Analytical data provided by {\sl YouTube}, we observe that the most popular videos that the students watch are the topics on Subspace, Linear Transformation and Finding the Jordan Canonical Form, where each video has been viewed for more than 600 times.

\subsection{Computer Algebra System {\sl SageMath}}

{\sl SageMath} is a CAS and it is also one of free open-source mathematics software. According to its official website, the mission of disseminating {\sl SageMath} is to create a viable free open source alternative to {\sl Magma}, {\sl Maple}, {\sl Mathematica} and {\sl Matlab}. See {\small \url{http://www.sagemath.org/}}. It is distributed under the terms of the GNU General Public License version 2$+$. It can be utilized to assist computations in Calculus, Linear Algebra, Combinatorics, Number Theory and Numerical Methods.

In this course of Introductory Linear Algebra, {\sl SageMath} is utilized as both in-class and out-of-class activities. In-class activities include demonstrations of solving Linear Algebra problems using {\sl SageMath} and problem-solving session for the students to try the CAS by themselves. Out-of-class activities include assignments for students to solve Linear Algebra problems using the CAS. For the former, the instructor and the students visit an online `single cell' version of {\sl SageMath} accessible at {\small \url{https://sagecell.sagemath.org/}} and perform a computation there.  The cloud computing platform is available from our own university's cell server and is accessible at {\small \url{sage.skku.edu}}. A screenshot of this cell server is shown on the right panel of Figure~\ref{sage}. Some alternative {\sl SageMath} cell servers available in South Korea amongst others are provided by Korea National Open University (KNOU), which is accessible at {\small \url{http://mathlab.knou.ac.kr:8080/}} and by the National Institute for Mathematical Sciences (NIMS), which can be accessed at {\small \url{https://sage.nims.re.kr/}}. The following subsection provides examples of flipped classroom activities by utilizing {\sl SageMath}.

\subsection{Flipped classroom activities}

Prior coming to the class, students are assigned to read a portion of the textbook on a relevant topic covering a particular week's material. They are also requested to watch video lectures on the relevant material, usually from {\sl YouTube}. The tasks are communicated via the LMS {\sl icampus}. After a brief summary of basic definition and theorem, a typical class time is dedicated to problem-solving session and workshop on the use of technology. Earlier in the semester, we demonstrate some computational examples using {\sl SageMath} to the students and they are also invited to attempt in solving problems using the software. The purpose is to help them to get more familiar with the CAS. We also select relatively simple examples with small matrix size so that the students can also confirm the computational results using hand computation. Later in the semester, the demonstration is conducted less frequent as the students have obtained a level of familiarity and they can figure out the syntax by themselves.
 
During the problem-solving session, the students can work as an individual or in a team that they form by themselves. On moving toward creating a collaborative learning environment, the latter is usually encouraged. We navigate around the classroom to monitor the students working on assigned problems. This occasion provides an opportunity not only for probing students' progress on the tasks and their understanding of Linear Algebra concepts but also in assisting them whenever they encounter difficulties. We observe that students are more open and less hesitant to ask questions during this session in comparison with the end-of-lecture session, where they are also given a chance to ask questions.
	
For a weaker group of students, we can personalize our approach by guiding them in solving problems with a slower step-by-step inquiry and explanation. They are also encouraged to practice more problems outside the problem-solving sessions and during out-of-class activities. For an academically stronger group of students and those who have a keen interest in mathematics, we provide them with more challenging and computationally demanding problems. We also encourage them to solve problems not only by hand computation but also using the CAS and ask them to confirm the result. In the case of some discrepancy, we inquire them whether there exists potentially some computational errors or the results are simply a matter of different representations of the solution. This may lead to further discussion in Linear Algebra, including conceptual one. Hence, in the flipped classroom, we tailor our pedagogical approach in satisfying the need of students with a diverse academic background. In either case, our role as instructors is to `guide on the side' rather than to be a `sage on the stage'~\cite{king1993sage}.

The following provides some examples of in-class and out-of-class problem-solving activities with the use of CAS {\sl SageMath}.
We start with the smaller size of matrices and then we gradually increase to a larger size of matrices.
We assign the students to find a determinant of $2 \times 2$ or $3 \times 3$ matrices, first using pen and paper, and then using the CAS.
For some special cases, a larger size of matrices can also be found by hand calculation, but CAS {\sl SageMath} gives a faster result, even for any size of matrices. Consider the following example of calculating the determinants of $3 \times 3$ and $5 \times 5$ matrices $A$ and $B$, respectively
\begin{equation*}
A = \left[
\begin{array}{crr}
1 &  5 &  0 \\
2 &  4 & -1 \\
0 & -2 &  0
\end{array}
\right] \qquad \qquad 
B = \left[
\begin{array}{crrrr}
3 & -7 &  8 &  9 & -6 \\
0 &  2 & -5 &  7 &  3 \\
0 &  0 &  1 &  5 &  0 \\
0 &  0 &  2 &  4 & -1 \\
0 &  0 &  0 & -2 &  0
\end{array}
\right].
\end{equation*}
The fastest way to calculate the determinant of matrix $A$ is by using the cofactor expansion across either the third row or the third column. To calculate the determinant of matrix $B$, one must recognize that the block of the final three rows and three columns of matrix $B$ contains matrix $A$. Hence, one can calculate the determinant of matrix $B$ by using the cofactor expansion along the first column, and then expanding along the second column and after that multiplying with the determinant of matrix $A$. Even though the students have calculated the determinant of matrix $A$ before, some of them fail to recognize the matrix block when calculating the determinant of matrix $B$. Using CAS {\sl SageMath}, the computation of these determinants only takes a few seconds. The following provides commands for the determinant computation.
{\small
\begin{verbatim}
A = matrix([[1,5,0], [2,4,-1], [0,-2,0]])
print "det(A) =", det(A)
B = matrix([[3,-7,8,9,-6],[0,2,-5,7,3],[0,0,1,5,0],[0,0,2,4,-1],[0,0,0,-2,0]])
print "det(B) =", det(B)
\end{verbatim}
}

One example of out-of-class activities is an assignment for the students to adopt CAS {\sl SageMath} or any other CAS they are familiar to work with in finding a dominant eigenvalue and the eigenvector corresponding to that particular dominant eigenvalue of the following $2 \times 2$ matrix $A$ with an initial condition $\mathbf{x}_0$, respectively
\begin{equation*}
A = \left[ \begin{array}{rr}
-7 & -12 \\ 8 & 13 
\end{array}\right] \qquad \qquad 
\mathbf{x}_0 = \left[\begin{array}{c}
1 \\ 0
\end{array} \right].
\end{equation*}
The students need to solve this problem numerically by implementing the `power method' algorithm and they were required to provide a minimum of ten number of iterations. A dominant eigenvalue of a matrix is the largest eigenvalue in the absolute value sense and the dominant eigenvector is the eigenvector corresponding to the dominant value. An example similar to this problem can be found in the free e-book, including the {\sl SageMath} code and the corresponding computational result. A common mistake is that the students calculate eigenvalues and the corresponding eigenvectors of the matrix analytically, instead of numerically, using {\sl SageMath} or any other CAS. Another common mistake that was observed is that although the students attempt to use power method, they did not provide a minimum ten number of iterations. The following provides a {\sl SageMath} code for finding the dominant eigenvalue and the eigenvector corresponding to this dominant eigenvalue for the example above.
{\small
\begin{verbatim}
from numpy import argmax,argmin
A=matrix([[-7,-12],[8,13]])
x0=vector([1.0,0.0]) # Initial guess of eigenvector
maxit=20 # Maximum number of iterates
dig=8 # number of decimal places to be shown is dig-1
tol=1e-7 # Tolerance limit for difference of two consecutive eigenvectors
err=1 # Initialization of tolerance
i=0
while(i<=maxit and err>=tol):
  y0=A*x0
  ymod=y0.apply_map(abs)
  imax=argmax(ymod)
  c1=y0[imax]
  x1=y0/c1
  err=norm(x0-x1)
  i=i+1
  x0=x1
  print "Iteration Number:", i-1
  print "y"+str(i-1)+" =",y0.n(digits=dig), 
        "c"+str(i-1)+" =", c1.n(digits=dig), 
        "x"+str(i)+" =",x0.n(digits=dig)
  print
\end{verbatim}}
{\setlength{\parindent}{1pt}
The computational result is displayed as follows.}
{\small
\begin{verbatim}
Iteration Number: 0
y0 = (-7.0000000, 8.0000000) c0 = 8.0000000 x1 = (-0.87500000, 1.0000000)

Iteration Number: 1
y1 = (-5.8750000, 6.0000000) c1 = 6.0000000 x2 = (-0.97916667, 1.0000000)

Iteration Number: 2
y2 = (-5.1458333, 5.1666667) c2 = 5.1666667 x3 = (-0.99596774, 1.0000000)

Iteration Number: 3
y3 = (-5.0282258, 5.0322581) c3 = 5.0322581 x4 = (-0.99919872, 1.0000000)

Iteration Number: 4
y4 = (-5.0056090, 5.0064103) c4 = 5.0064103 x5 = (-0.99983995, 1.0000000)

Iteration Number: 5
y5 = (-5.0011204, 5.0012804) c5 = 5.0012804 x6 = (-0.99996800, 1.0000000)

Iteration Number: 6
y6 = (-5.0002240, 5.0002560) c6 = 5.0002560 x7 = (-0.99999360, 1.0000000)

Iteration Number: 7
y7 = (-5.0000448, 5.0000512) c7 = 5.0000512 x8 = (-0.99999872, 1.0000000)

Iteration Number: 8
y8 = (-5.0000090, 5.0000102) c8 = 5.0000102 x9 = (-0.99999974, 1.0000000)

Iteration Number: 9
y9 = (-5.0000018, 5.0000020) c9 = 5.0000020 x10 = (-0.99999995, 1.0000000)

Iteration Number: 10
y10 = (-5.0000004, 5.0000004) c10 = 5.0000004 x11 = (-0.99999999, 1.0000000)
\end{verbatim}
}From this result, we observe that the dominant eigenvalue is 5 and the eigenvector corresponding to the dominant eigenvalue is $(-1, 1)$.

\section{Students' Feedback and Reflection} \label{feedback}

After implementing flipped classroom approach in our Linear Algebra classes, we are interested in investigating what our students think and feel about this pedagogy. In this section, we discuss the feedback from the students so that we can reflect on their concerns in order to improve for a better teaching and learning in a flipped classroom environment. 

\subsection{Students' feedback}

The feedback was collected in the middle of the semester of Spring 2016, on Week~7, one week prior to the midterm exam period. The students were asked on their opinion of the flipped classroom we have implemented so far. There also exists another feedback collected at nearly the end of the semester before the final exam period, around Weeks~14 and~15. A number of selected comments from the students is given as follows. 
{\small \begin{quoting}[leftmargin=0.25cm]
\textmd{``I hate a flipped class \dots. I want a traditional class.''} \\
\textrm{``No flipped class, please. Lecture in the classroom is more useful.''}\\
\textmd{``A flipped class is not good. I feel that the flipped class does not have a vitality.''}
\end{quoting}}
From these comments, we observe that some students do not really favor the flipped classroom approach and express negative reactions toward it. They prefer the traditional teaching style where an instructor is lecturing throughout the entire class session and the students are sitting down listening and taking notes. This is understandable since the traditional class has become a norm ever since they attended school from their childhood. Besides, they also think that  it is more convenient and much easier on their part by simply arriving at the classroom without any prior preparation. When they enroll in other classes which are not delivered as the flipped classroom, they adopt a similar habit of coming to the class unprepared. After all, enrolling in a flipped classroom course requires plenty of effort. The students do not only need to prepare prior coming to the class but they are also required to participate actively in the class. Indeed, before we are able to flip the class successfully, we also need to flip students' attitude and mindset in the first place. This is one of the challenges in implementing flipped classroom approach successfully. 

In connection with the video recordings, the following are some comments that the students have written.
{\small \begin{quoting}[leftmargin=0.25cm]
\textrm{``It is hard to understand what the video says because it is too short, so I have to study alone. It feels that I don't take any class but just to studying alone.''} \\
\textmd{``The videos are just reading presentation slide. Please add more explanation.''} \\
\textrm{`` \dots it was not good that the lecture was given in the video. In the video, the professor just reading the presentation slide.''} \\
\textmd{``I think just studying with recorded video is hard, so your teaching during the class is needful.''}\\
\textrm{``I need more `detailed explanation' in the recorded video. I could not understand the content fully unless I watch another video on \textsl{YouTube}, etc.''} \\
\textmd{``Video lectures are not helpful. So, it is hard to understand the concept of Linear Algebra.''}
\end{quoting}}
From these comments, we discover that a lot of things need to be improved when it comes to producing recorded video lectures. In a course like Introductory to Linear Algebra, we need to inform the students of many important definitions, terminologies, lemmas, propositions and theorems and we thought that informing this information through the video will be a good idea. Therefore, we attempt to read these definitions and theorems so that we do not need to repeat telling the same thing during the class time. However, students interpreted this action differently. They informed us that simply reading the texts is not a good thing to do. It may inform them, but does not necessarily they understand the presented material. Therefore, a more detailed explanation is needed and some examples accompanying the presentation will do good.

\subsection{Reflection} 

Implementing the flipped classroom approach certainly is not an easy task. The students need to be introduced clearly what to be expected from them at the beginning of the semester and also throughout the period repetitively. We need to convince the students that their efforts are worthed when they were asked to spend more time for class preparation outside the regular class contact hours. All of us are fully aware that our students are busy too with other classes and extracurricular activities. Therefore, it is not always easy for them to do preparation, even just merely watching the videos, on a regular basis for the entire period of the semester.

On a technical side, we realized that producing high-quality videos is not an easy task and teaching using video can be pretty hard as well. We did not receive any official training on how to produce good quality teaching video recordings. It takes several times of trials and errors. We also need to ensure that the duration of the video recordings is not too short and not too long either. In fact, it has been investigated that videos with a shorter duration are more engaging to the students~\cite{guo2014video}. For many of the students, particularly those who are used to with traditional style of teaching, receiving information for the first time through reading a textbook or watching the video recordings can be a challenge. They might have difficulty in grasping the presented materials, particularly if English is not their first language. Hence, more issues need to be addressed on implementing flipped classroom pedagogy in the context of English-medium instruction environment. We also observed that providing examples in our videos can assist the students to make a better connection between conceptual aspects and practical methods. This, in turn, will help them in solving Linear Algebra problems. 

\section{Conclusion} \label{conclusion}

In this article, we have discussed our experience in flipping the classroom of an Introductory Linear Algebra course at SKKU. A new feature of our teaching includes developing a new, free electronic textbook that can be accessible to everyone in the world and disseminating online lecture notes on a single web page that can be used by both students and instructors alike. We have embedded the use of technology in our teaching, in particular exploring the power of open access and free CAS {\sl SageMath} for computational purposes. For communication between instructors and students outside the classroom, we utilize an intranet LMS {\sl icampus}. This platform is also used to provide additional teaching materials, to send announcements and to foster discussion in an online forum.

We have also created recorded video lectures, either in a screencast format or in the form of real live classroom teaching. These videos are accessible from the video-sharing platform {\sl YouTube} and the students can watch them primarily not only before they come to the class sessions, but also after the class is over if they wish to review the lectures again at their own pace. The students are expected to gather basic information on definitions, properties, theorems, methods and techniques in Linear Algebra prior the class time, either by reading assignments, listening to podcasts or watching recorded video lectures.  We dedicated class time to activities that constitute active learning, including problem-solving, question and answer, discussion and providing feedback. We observe that by navigating around the classroom, we do not only discover students' hurdles in understanding concepts and techniques but also tailor our pedagogy in assisting students with different academic ability and in various stages of learning to learn Linear Algebra successfully.

Students' feedback and opinion on their experience in flipped classroom environment indicate a less favorable reaction. Some expressed that they prefer the traditional classroom instead of the flipped classroom. This can be explained by either they are used to with enrolling courses delivered in the traditional lecture style format or the way we design and implement the flipped classroom does not satisfy students' learning styles. For the latter, the students pinpoint the quality of the video recordings that we have created. The mentioned that the videos are not really interesting if they are only reading through definitions and theorems instead of working through some examples. Despite this technical aspects, some students tried their best to participate actively in the problem-solving sessions and we have an opportunity to communicate with them. Overall, we have gleaned valuable lessons in flipping the class so that we can improve and implement a better strategy in the flipped classroom pedagogy for future courses we are going to teach.

\subsection*{Acknowledgement}
{\small The authors would like to acknowledge SKKU's Center of Teaching and Learning (CTL) for the financial support in enabling us to implement the flipped classroom pedagogy in our Linear Algebra courses. SGL acknowledged the support from Basic Science Research Program through the National Research Foundation of Korea (NRF) funded by the Korean Ministry of Education (2017R1D1A1B03035865). We express gratitude to Dr. Sang-Eun Lee (CTL, SKKU),  Dr. Rob Layahe (Department of Physics, University College, SKKU), Dr. Lois Simon (Department of Mathematics, University College, SKKU), Professor Marco Pollanen (Trent University, Canada) for advice, feedback, fruitful discussions on various practical aspects involved in implementing flipped classroom pedagogy. We would like to thank our students for their honest feedback and willingness in trying a new learning environment. Finally, we thank anonymous reviewers for their constructive comments and suggestions which have improved the quality of this paper significantly.  \par}

{\small
\bibliography{18eJMTeBook} 
\bibliographystyle{plain}
}

\end{document}